\documentclass[a4paper,11pt]{article}

\usepackage[pdfstartview = {FitH}, bookmarksopen]{hyperref} 

\usepackage{amsfonts} 

\usepackage[latin1]{inputenc}  

\usepackage{amssymb}   					
\usepackage{mathrsfs}                    
\usepackage{amsmath}                    
\usepackage{amsthm} 
\usepackage{xcolor}                    
\usepackage[arrow, matrix, curve]{xy} 

\usepackage{color}

\usepackage{pdfpages}
\usepackage{multirow}   

\usepackage{rotating}
\usepackage{rotfloat}
\usepackage{lscape} 

\usepackage{tikz}
\usetikzlibrary{shapes,arrows,calc}
\tikzstyle{decision} = [chamfered rectangle, draw, text centered, chamfered rectangle xsep=2em]
\tikzstyle{block}    = [rectangle, draw, text centered]
\tikzstyle{schrift}  = [rectangle, text centered]
\tikzstyle{line}     = [draw, -latex']  
\tikzstyle{symbol}   = [rectangle, draw, text centered, minimum height= 2em, minimum width = 2em]
\tikzstyle{syndrom}  = [rectangle, draw, text centered, minimum height= 2em, minimum width = 4em]

\theoremstyle{plain}
\newtheorem{theorem}{Theorem}[section] 
\newtheorem{proposition}[theorem]{Proposition}

\newtheorem{example}[theorem]{Example}
\newtheorem{remark}[theorem]{Remark}

\theoremstyle{definition}
\newtheorem{definition}[theorem]{Definition}

\title{Cryptographic protocols based on Nielsen transformations}
\author{Anja I. S. Moldenhauer and Gerhard Rosenberger}
\date{}  

\begin{document}

\maketitle

\begin{abstract}
Based on a combinatorial distribution of shares we present in this paper  secret sharing schemes and cryptosystems using Nielsen transformations. 
\end{abstract}




\section{Introduction}\label{introPan}

We first describe secret sharing protocols and combinatorial distributions of shares. After this introductory definitions we start with  a secret sharing scheme using directly the combinatorial distribution of shares. Based on this we present two schemes in which we apply regular Nielsen transformations in connections with faithful representations of free groups  and the Nielsen reduction theory. In the last sections we modify the secret sharing schemes to a  private key cryptosystem and finally  Nielsen transformations are used for a public key cryptosystem which is inspired by the ElGamal cryptosystem. The new cryptographic protocols are in part in the dissertation from A. Moldenhauer \cite{M} under her supervisor G. Rosenberger at the University of Hamburg. \\

A $(n,t)$-secret sharing protocol, with $n,t \in \mathbb N$ and $t \le  n$, is a method to distribute a secret among a group of $n$ participants in such a way that it can be recovered
only if at least $t$ of them combine their shares. Hence any group of $t-1$ or fewer participants cannot calculate the secret. The number $t$ is called threshold.
The person who distributes the shares is called the dealer.

D. Panagopoulos presents in his paper \cite{P} a $(n,t)$-secret sharing scheme using group presentations with solvable word problem. Here we use  combinatorial distributions of the shares similar to those introduced in the paper of  D.~Panagopoulos: 

To distribute the shares in a $(n,t)$-secret sharing scheme the dealer does the following steps: 
 
\begin{enumerate}
\item Calculate  $m=\binom{n}{t-1}$, the number of all elements, for example $\{a_1,a_2,\ldots,a_m\}$, the participants need to know for the  reconstruction of the  secret.
\item Let $A_1,A_2,\ldots,A_m$ be an enumeration of the subsets of $\{1,2,\ldots,n\}$ with $t-1$ elements. Define $n$ subsets $R_1,R_2,\ldots,R_n$ of the set $\{a_1,a_2,\ldots,a_m\}$ with the property
      \[
      a_j \in R_i  \qquad \Longleftrightarrow \qquad i \not\in A_j \qquad \text{for } j=1, 2, \ldots, m \text{ and } i=1, 2, \ldots, n.
      \]
\item The dealer distributes to each of the $n$ participants one of the sets $R_1,R_2,\ldots,R_n$.
\end{enumerate}

The new protocols in this paper are based on Nielsen transformations, which are the basis of a linear technique to study free groups and  general infinite groups. We now review some basic definitions concerning regular Nielsen transformations and Nielsen reduced sets (see \cite{CRR} or \cite{LS}).\\

Let $F$ be a free group on the free generating set $X:=\{x_1,x_2,\ldots \}$ and let $U:=\{u_1,u_2,\ldots \} \subset F$.

\begin{definition}\label{Nreduceddef}
An \textbf{elementary Nielsen transformation} on \newline  $U=\{u_1,u_2,\ldots \}$ is one of the following transformations
\begin{enumerate}
\item[(T1)] replace some $u_i$ by $u_i^{-1}$;
\item[(T2)] replace some $u_i$ by $u_iu_j$ where $j \ne i$;
\item[(T3)] delete some $u_i$ where $u_i=1$.
\end{enumerate}
In all three cases the $u_k$ for $i \ne k$ are not changed.
A (finite) product of elementary Nielsen transformations is called a \textbf{Nielsen transformation}.
A Nielsen transformation is called \textbf{regular} if it is a finite product of the transformations $(T1)$ and $(T2)$, otherwise it is called \textbf{singular}.
The set $U$ is called \textbf{Nielsen-equivalent} to the set $V$, if there is a regular Nielsen transformation from $U$ to $V$.
\end{definition}

\begin{definition}
Consider elements $v_1,v_2,v_3$ of the form $u_i^{\pm 1}$, call $U$ \textbf{Nielsen reduced} if for all such triples the following conditions hold:
\begin{enumerate}
\item[(N0)] $v_1 \ne 1$;
\item[(N1)] $v_1v_2 \ne 1$ implies $|v_1v_2|\geq |v_1|, |v_2|$; 
\item[(N2)] $v_1v_2\ne 1$ and $v_2v_3 \ne 1$ implies $|v_1v_2v_3|>|v_1| - |v_2| + |v_3|$.
\end{enumerate}
Here $|\cdot|$ denotes the free length in $F$.
\end{definition}

\begin{proposition}\label{propNred}
If $U=\{u_1,u_2,\ldots,u_n\}$ is finite, then $U$ can be carried by a Nielsen transformation into some $V$ such that $V$ is Nielsen reduced.
\end{proposition}

For a proof see \cite[Theorem 2.3]{CRR} or \cite[Proposition 2.2]{LS}.

For the secret sharing scheme based on Nielsen transformations we will only use regular Nielsen transformations.
We agree on some notations.

We write $(T1)_i$ if we replace $u_i$ by $u_i^{-1}$ and we write $(T2)_{ij}$ if we replace $u_i$ by $u_iu_j$. If we want to apply $t$-times one after the other the same Nielsen transformation $(T2)$ we write $[(T2)_{ij}]^t$ and hence replace $u_i$ by $u_iu_j^t$. In all cases the $u_k$ for $i\ne k$ are not changed.


\section{A combinatorial secret sharing scheme}\label{combsecretscharing}

Now we present a $(n,t)$-secret sharing scheme, whereby the secret is the sum 
of multiplicative inverses of elements in the natural numbers.
For the distribution of the shares the dealer uses the method 
of D.~Panago\-poulos described in Section \ref{introPan}.

The numbers $n$ and $t$ are given, whereby $n$ is the number of participants and $t$ is the threshold.

\begin{enumerate}
\item The dealer first calculates the number $m=\binom{n}{t-1}$. 
\item He chooses $m$ elements $a_1, a_2, \ldots, a_m \in \mathbb N$.  From these elements he constructs analogously as in Section \ref{introPan} the sets $R_1, R_2, \ldots ,R_n$. 
The secret $S$ is the sum
\begin{align*}
S:= \sum_{i=1}^{m} \frac{1}{a_i}  \in \mathbb Q^{+}.
\end{align*}
\item Each participant $P_i$ gets one share $R_i$, $1 \leq i \leq n$.
\end{enumerate}


If $t$ of the $n$ participants come together they can reconstruct the secret while they first combine their $t$ private sets $R_i$ and get by construction the set $\tilde{R}~=~\{a_1, a_2, \ldots, a_m\}$.
The secret is the sum of the inverse elements in the set $\tilde{R}$, that is
\begin{align*}
 S =  \sum_{i=1}^{m} \frac{1}{a_i}.
\end{align*}

If the dealer needs a special secret $\tilde{S} \in \mathbb Q$  he gives every participant one more element $x \in \mathbb Q$ in each $R_i$, with
\begin{align*}
x:= \frac{\tilde{S}}{S}.
\end{align*}
The participants get $\tilde{S}$ by multiplying the reconstructed secret $S$ with $x$.


Each element $a_j$ is exactly contained in $n-(t-1)$ subsets. Hence for each $j=1,2,\ldots,m$ the element $a_j$ is not contained in $t-1$ subsets from $\{R_1,R_2,\ldots,R_n\}$. As a consequence, $a_j$ is in each union of $t$ subsets. Otherwise, if just $t-1$ arbitrary sets from $\{R_1,R_2, \ldots, R_n\}$ are combined, there exist  a $j$ so that the element $a_j$ is not included in the union of this sets.

If just one element $a_j$ is absent, the participants do not get the correct sum~$S$, 
and hence cannot compute the correct secret.

\begin{example} 

We perform the steps for a $(4,3)$-secret sharing scheme.
It is $n=4$ and $t=3$.

The dealer follows the steps:

\begin{enumerate}
\item He first calculates $m=\binom{n}{t-1}=\binom{4}{2}= 6$.
\item The dealer chooses the numbers $a_1:=2, a_2:=1,a_3:=2,a_4:=8$, \newline $a_5:=4$ and $a_6:=2$. The secret is 
 \begin{align*}
    S:= \sum_{i=1}^{m} \frac{1}{a_i} = \frac{23}{8}.
 \end{align*}
\begin{enumerate}
\item The six subsets with size 2 of the set $\left\{1, 2, 3, 4 \right\}$
are
\begin{alignat*}{4}
      A_1 &= \left\{1, 2\right\}, \qquad A_2 &= \left\{1, 3\right\}, \qquad A_3 &= \left\{1, 4\right\},\\
      A_4 &= \left\{2, 3\right\}, \qquad A_5 &= \left\{2, 4\right\}, \qquad A_6 &= \left\{3, 4\right\}.
      \end{alignat*}
     With help of the $A_i$ the dealer gets the sets $R_1, R_2, R_3$ and $R_4$, which contain elements from $\left\{a_1, \ldots, a_6\right\}$. He puts the element $a_{j}$ for which $i$ is not contained in the set $A_{j}$ for  
      $i=1,\ldots,4$  
      and $j=1,\ldots,6$, into the set $R_{i}$:
      \begin{align*}
      1 \not\in A_4, A_5, A_6  \Longrightarrow  R_1=\left\{a_4, a_5, a_6\right\},\\
      2 \not\in A_2, A_3, A_6  \Longrightarrow  R_2=\left\{a_2, a_3, a_6\right\},\\
      3 \not\in A_1, A_3, A_5  \Longrightarrow  R_3=\left\{a_1, a_3, a_5\right\},\\
      4 \not\in A_1, A_2, A_4  \Longrightarrow  R_4=\left\{a_1, a_2, a_4\right\}.
      \end{align*}
\end{enumerate}
\item The dealer distributes the set $R_i$  to the participant $T_i$, for $i=1,\ldots,4$.
\end{enumerate}

If three of the four participants come together, they can calculate the secret~$S$. For example the participants $T_1, T_2$ and $T_3$ have the set
\begin{align*} \tilde{R}:=& R_1 \cup R_2 \cup R_3\\
=&\left\{a_4, a_5, a_6\right\} \cup \left\{a_2, a_3, a_6\right\} \cup \left\{a_1, a_3, a_5\right\}\\
                =& \left\{a_1, a_2,a_3, a_4, a_5, a_6\right\},
\end{align*}

and hence get the secret 
 \begin{align*}
    S= \sum_{i=1}^{6} \frac{1}{a_i} = \frac{23}{8} \text{\qquad  with } a_i \in \tilde{R}.
 \end{align*}

\end{example}


\section{A secret sharing scheme using a regular Nielsen transformation}\label{NSec}

In this section we describe a $(n,t)$-secret sharing scheme which extends and improves the ideas in Section \ref{combsecretscharing} by using Nielsen transformations. We consider free groups as abstract groups but also as subgroups of the special linear group of all $2 \times 2$ matrices over $\mathbb Q$, that is,
\begin{align*}
SL (2,\mathbb Q)=\left\{ \begin{pmatrix}
 a & b\\
 c & d
\end{pmatrix} ; a,b,c,d \in \mathbb Q \text{ and } ad-bc=1 \right\}.
\end{align*}

We use the special linear group over the rational numbers because these numbers can be stored and computed more efficiently on a computer than irrational numbers.

Let $F$ be a free group in $SL(2,\mathbb Q)$ of rank $m:=\binom{n}{t-1}$. The dealer wants to distribute the shares for the participants as described in Section \ref{introPan}. The shares will be subsets of a free generating set of the group $F$.


\textbf{Steps for the Dealer:}
The numbers $n$ and $t$ are given, whereby $n$ is the number of participants and $t$ is the threshold.
We have $m:=\binom{n}{t-1}$.
\begin{enumerate}
\item The dealer chooses an abstract free generating set $X$ for the free group $F$ of rank $m$, that is
	$$F= \langle X; \phantom{R} \rangle \text{\qquad  with }  X:=\{ x_1, x_2,\ldots, x_m\}.$$
	He also needs an explicit free generating set $M$, that is
	$$F= \langle M; \phantom{R} \rangle \text{\qquad  with }  M:=\{ M_1, M_2,\ldots, M_m \}$$
	and $M_i \in 		SL(2,\mathbb Q)$.
\item With the known matrices in the set $M$ he computes the secret 
	\begin{align*}
	S:= \sum_{j=1}^{m} \frac{1}{|a_j|} \in \mathbb Q^{+} \qquad \text{ with } a_j:=tr(M_j) \in \mathbb Q,
	\end{align*}
	$tr(M_j)$ is the trace for the matrix $M_i:=\begin{pmatrix}
	a & b\\
	c & d
	\end{pmatrix} \in SL(2,\mathbb Q)$, that is, $tr(M_i):=a+d$. 
	If the dealer needs a special secret he can act as in Section \ref{combsecretscharing} described.
\item The dealer constructs the shares for the participants in the following way:
	\begin{enumerate}
		\item He first applies a regular Nielsen transformation simultaneously for both sets $X$ and $M$ to
			get Nielsen-equivalent sets $U$ and $N$ to $X$ and $M$, respectively (see Figure \ref{SimNilsfor}).

			\begin{figure}[H]
			\begin{center}			
			\hspace{1cm} $X:=\{ x_1, x_2,\ldots, x_m\}$ \hspace{3cm} $M:=\{ M_1, M_2,\ldots, M_m \}$\\
			\begin{tikzpicture}
			\draw[->] (1,1)--(1,-1);
			\node (s) at (-0.5,0) {regular Nielsen};
			\node (t) at (-0.5,-0.5) {transformation};
			\draw[->] (8,1)--(8,-1);
			\node (u) at (6.5,0) {regular Nielsen};
			\node (v) at (6.5,-0.5) {transformation};
			\end{tikzpicture}\\
			\hspace{1cm} $U:=\{ u_1, u_2,\ldots, u_m\}$  
			\hspace{3cm} $N:=\{ N_1, N_2,\ldots, N_m \}$
			\caption{Simultaneously regular Nielsen transformation}\label{SimNilsfor}
			\end{center}
			\end{figure}
		
			The elements $u_i$ are words in $X$ and the elements $N_i$ are words in~$M$. 
			Hence we have $N_i~\in~SL(2,\mathbb Q)$.
			
		\item The dealer now uses the method of D.~Panagopoulos to split $U$ and $N$ and to get the shares 
		$(R_i,S_j)$ for the participants with $R_i \subset U$ and $S_j \subset N$.
	\end{enumerate}
	\item The dealer distributes the shares.
\end{enumerate}


If $t$ of the $n$ participants combine their parts  they obtain the sets $U$ and $N$. 
The secret can be recovered as follows:
\begin{enumerate}
\item The participants apply regular Nielsen transformations in a Nielsen reduction manner for $U$ and step by step simultaneously for $N$. By Proposition \ref{propNred} they get Nielsen reduced sets $X^{\pm}= \{x_1^{\epsilon_1}, x_2^{\epsilon_2}, \ldots, x_m^{\epsilon_m} \}$ and $M^{\pm}=\{M_1^{\delta_1},M_2^{\delta_2},\ldots,M_m^{\delta_m}\}$ with $\epsilon_i, \delta_i \in \{+1,-1\}$, see Figure \ref{SimNilsback}.

			\begin{figure}[H]
			\begin{center}
			\hspace{1cm} $U:=\{ u_1, u_2,\ldots, u_m\}$  \hspace{3cm} $N:=\{ N_1, N_2,\ldots, N_m \}$ \\
			\begin{tikzpicture}
			\draw[->] (1,1)--(1,-1);
			\node (s) at (-0.5,0) {regular Nielsen};
			\node (t) at (-0.5,-0.5) {transformations};
			\draw[->] (8,1)--(8,-1);
			\node (u) at (6.5,0) {regular Nielsen};
			\node (v) at (6.5,-0.5) {transformations};
			\end{tikzpicture}\\
			
			\hspace{1cm} $X^{\pm}= \{x_1^{\epsilon_1}, x_2^{\epsilon_2}, \ldots, x_m^{\epsilon_m} \}$ \hspace{3cm} $M^{\pm}=\{M_1^{\delta_1},M_2^{\delta_2},\ldots,M_m^{\delta_m}\}$
			\caption{Simultaneously regular Nielsen transformations}\label{SimNilsback}
			\end{center}
			\end{figure}
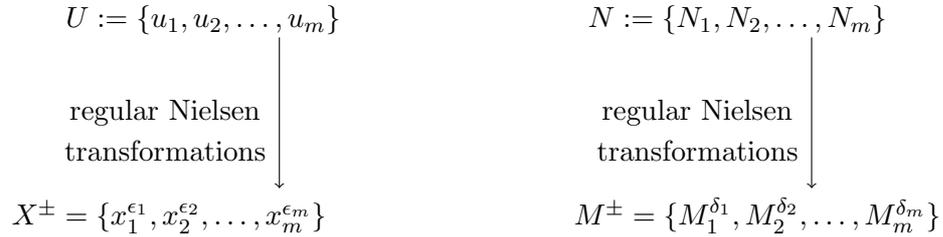
\item With the knowledge of the set $M^{\pm}$ it is easy to reconstruct the secret
\begin{align*}
	S=  \sum_{j=1}^{m} \frac{1}{|a_j|} \in \mathbb Q^{+} \qquad \text{ with } tr(M_j)=a_j \in \mathbb Q.
	\end{align*}
	Recall that $tr(M_i^{\delta_i})=tr(M_i)$ for $i=1,\ldots,m$.
\end{enumerate}


Less than $t$ participants can neither  get the whole set $U$, which is Nielsen-equivalent to $X$, nor the set $N$, which is Nielsen-equivalent to $M$.\\
For the calculation of the secret, the participants need the set $M$, because the secret depends on the traces of the matrices $M_i \in M$.
The participants need both sets $U$ and $N$. If they just have one set $U$ or $N$ they cannot get information about the set $M$.\\
 If the set $U$ is known, it is only known which Nielsen transformation should be done to get the Nielsen-equivalent set $X$, but it is unknown on which matrices they should be done simultaneously.\\
If only the set $N$ is known, then the matrices in $SL(2,\mathbb Q)$ are known, but nobody knows which Nielsen transformation should be done on $N$ to get the set $M$. It is also unknown how many Nielsen transformations were used.


In the book \cite{L} of J. Lehner  on page 247 a method is given to explicitly obtain a 
free generating set $M$ for a free group $F$ on the abstract generating set 
$X~:=~\{ x_1, x_2,\ldots, x_m\}$:

\begin{example}\label{examplmatrices}
Let $F$ be a free group with countably many free generators $x_1,x_2,\ldots$. Corresponding to $x_j$ define the matrix 
\begin{align*}
M_j=\begin{pmatrix}
 -r_j & -1+r_j^2\\
 1 & -r_j
\end{pmatrix}
\end{align*}
with $r_j \in \mathbb Q$ such  that the following inequalities hold:
\begin{align}\label{inequalitis}
r_{j+1}-r_j \geq 3 \qquad \text{ and } \qquad r_1\geq 2.
\end{align}
The group $G$ generated by $\{M_1,M_2,\ldots\}$ is isomorphic to $F$ (see \cite{L}).
\end{example}

We now present an example for this secret sharing scheme.

\begin{example}
We perform the steps for  a $(3,2)$-secret sharing scheme with the help of the computer program Maple 16.
It is $n=3$, $t=2$ and hence $m~=~\binom{3}{1}=3$.

First the Dealer generates the shares for the participants.
\begin{enumerate}
\item The dealer chooses an abstract presentation for the free group $F$ of rank~$3$
	$$F= \langle X; \phantom{R} \rangle \text{\qquad  with }  X:=\{ x_1, x_2, x_3\}.$$
	He takes an explicit presentation
	$$F= \langle M; \phantom{R} \rangle \text{\qquad  with }  M:=\{ M_1, M_2,M_3 \},$$
	 $M_i \in SL(2,\mathbb Q)$ as above.  We first mention that the inequalities \eqref{inequalitis} hold for 
	\begin{align*}
	r_1= \frac{7}{2}, \quad r_2= \frac{15}{2}, \quad r_3= 11
	\end{align*}
	 and hence the set of the matrices
	\begin{align*}
		M_1&=\begin{pmatrix}
		-\frac{7}{2} & -1+\left( \frac{7}{2}\right)^2\\
		1 & -\frac{7}{2}
		\end{pmatrix}=\begin{pmatrix}
		-\frac{7}{2} & \frac{45}{4}\\
		1 & -\frac{7}{2}
		\end{pmatrix},\\
		M_2&=\begin{pmatrix}
		-\frac{15}{2} & -1+\left(\frac{15}{2}\right)^2\\
		1 & -\frac{15}{2}
		\end{pmatrix}=\begin{pmatrix}
		-\frac{15}{2} & \frac{221}{4}\\
		1 & -\frac{15}{2}
		\end{pmatrix},\\
		M_3&=\begin{pmatrix}
		-11 & -1+11^2\\
		1 & -11
		\end{pmatrix}=\begin{pmatrix}
		-11 & 120\\
		1 & -11
		\end{pmatrix}
	\end{align*}	
	is a free generating set for a free group of rank $3$.
\item We have
	\begin{align*}
	a_1:=tr(M_1)=-7, \quad a_2:=tr(M_2)=-15, \quad a_3:=tr(M_3)=-22,
	\end{align*}
	and hence the secret is
	\begin{align*}
	S:=  \sum_{j=1}^{3} \frac{1}{|a_j|}  = \frac{589}{2310}.
	\end{align*}
	\item Construction of the shares for the participants:
	\begin{enumerate}
	
\item First the dealer applies  regular Nielsen transformations (NTs) simultaneously for both sets $X$ and $M$ to
			get Nielsen-equivalent sets $U$ and $N$ to $X$ or $M$, respectively. 
			These transformations are shown in the Table \ref{vorgehen}.

	\begin{table}[H]
	\begin{small}
	\begin{center}
	\setlength{\tabcolsep}{.3\tabcolsep}
	\caption{Nielsen transformations (NTs) of the dealer}\label{vorgehen}

	\begin{tabular}{|l | l | l|}
	\hline
	& & \\
	NTs & theoretical set A  & \hspace{1cm} explicit set M \\
	& & \\
	\hline
	\hline
		& & \\
		& $\{x_1,x_2,x_3\}$ & $\Biggl\{\begin{pmatrix}
		-\frac{7}{2} & \frac{45}{4}\\
		1 & -\frac{7}{2}
		\end{pmatrix},
		
		\begin{pmatrix}
		-\frac{15}{2} & \frac{221}{4}\\
		1 & -\frac{15}{2}
		\end{pmatrix}$,
		
		$\begin{pmatrix}
		-11 & 120\\
		1 & -11
		\end{pmatrix} \Biggr\}$\\
	& & \\
	
	\hline
		& & \\
	$(T1)_2$ & $\{x_1,x^{-1}_2,x_3\}$ & $\Biggl\{\begin{pmatrix}
		-\frac{7}{2} & \frac{45}{4}\\
		1 & -\frac{7}{2}
		\end{pmatrix},
		
		\begin{pmatrix}
		-\frac{15}{2} & -\frac{221}{4}\\
		-1 & -\frac{15}{2}
		\end{pmatrix}$,
		
		$\begin{pmatrix}
		-11 & 120\\
		1 & -11
		\end{pmatrix} \Biggr\}$\\
	& & \\

	\hline
	& & \\
	$(T2)_{1.2}$ & $\{x_1x_2^{-1},x^{-1}_2,x_3\}$ & $\Biggl\{
	    \left(\begin{smallmatrix}
		15 & 109\\
		-4 & -29
		\end{smallmatrix}\right),
		
		\left(\begin{smallmatrix}
		-\frac{15}{2} & -\frac{221}{4}\\
		-1 & -\frac{15}{2}
		\end{smallmatrix}\right)$,
		
		$\left(\begin{smallmatrix}
		-11 & 120\\
		1 & -11
		\end{smallmatrix}\right) \Biggr\}$\\
	& & \\
	\hline	
	
		& & \\
	$\lbrack(T2)_{3.2}\rbrack^3$ & $\{x_1x^{-1}_2,x^{-1}_2,x_3x^{-3}_2\}$ & $\Biggl\{
	   \left( \begin{smallmatrix}
		15 & 109\\
		-4 & -29
		\end{smallmatrix}\right),
		
		\left( \begin{smallmatrix}
		-\frac{15}{2} & -\frac{221}{4}\\
		-1 & -\frac{15}{2}
		\end{smallmatrix}\right)$,
		
		$\left(\begin{smallmatrix}
		-8565 & -63664 \\
		799 & 5939
		\end{smallmatrix}\right) \Biggr\}$\\
	& & \\
	\hline

		& & \\
	$(T2)_{2.3}$ & $\{x_1x^{-1}_2,x^{-1}_2x_3x^{-3}_2,x_3x^{-3}_2\}$ & $\Biggl\{
	    \left( \begin{smallmatrix}
		15 & 109\\
		-4 & -29
		\end{smallmatrix}\right),
		
		\left(\begin{smallmatrix}
		\frac{80371}{4} & \frac{597401}{4}\\
		\frac{5145}{2} & \frac{38243}{2}
		\end{smallmatrix}\right)$,
		
		$\left(\begin{smallmatrix}
		-8565 & -63664 \\
		799 & 5939
		\end{smallmatrix}\right) \Biggr\}$\\
	& & \\
	\hline

		& & \\
	$(T1)_{1}$ & $\{x_2x^{-1}_1,x^{-1}_2x_3x^{-3}_2,x_3x^{-3}_2\}$ & $\Biggl\{
	    \left( \begin{smallmatrix}
		-29 & -109\\
		4 & 15
		\end{smallmatrix}\right),\!
		
		\left( \begin{smallmatrix}
		\frac{80371}{4} & \frac{597401}{4}\\
		\frac{5145}{2} & \frac{38243}{2}
		\end{smallmatrix}\right)$,\!
		
		$\left( \begin{smallmatrix}
		-8565 & -63664 \\
		799 & 5939
		\end{smallmatrix}\right)\! \Biggr\}$\\
	& & \\
	\hline
	
		& & \\
	$(T2)_{1.2}$ & $\{x_2x^{-1}_1x^{-1}_2x_3x^{-3}_2,$ &
	$\Biggl\{ \left( \begin{smallmatrix}
		-\frac{3452369}{4} & -\frac{25661603}{4}\\
		\frac{237917}{2} & \frac{1768447}{2}
		\end{smallmatrix}\right)$,\\
		 & $x^{-1}_2x_3x^{-3}_2,x_3x^{-3}_2\}$ & 
		$\left( \begin{smallmatrix}
		\frac{80371}{4} & \frac{597401}{4}\\
		\frac{5145}{2} & \frac{38243}{2}
		\end{smallmatrix}\right)$,
		
		$\left( \begin{smallmatrix}
		-8565 & -63664 \\
		799 & 5939
		\end{smallmatrix}\right) \Biggr\}$\\
	& & \\
\hline

	& & \\
	$(T1)_{3}$ & $\{x_2x^{-1}_1x^{-1}_2x_3x^{-3}_2,$ &
	$\Biggl\{ \left( \begin{smallmatrix}
		-\frac{3452369}{4} & -\frac{25661603}{4}\\
		\frac{237917}{2} & \frac{1768447}{2}
		\end{smallmatrix}\right)$,\\
		
		 & $x^{-1}_2x_3x^{-3}_2,x^3_2x^{-1}_3\}$ & 
		$\left( \begin{smallmatrix}
		\frac{80371}{4} & \frac{597401}{4}\\
		\frac{5145}{2} & \frac{38243}{2}
		\end{smallmatrix}\right)$,
		
		$\left( \begin{smallmatrix}
		5939 & 63664\\
		-799 & -8565
		\end{smallmatrix}\right) \Biggr\}$\\
	& & \\

\hline	
	& & \\
	$(T2)_{3.2}$ & $\{x_2x^{-1}_1x^{-1}_2x_3x^{-3}_2,$ &
	$\Biggl\{ \left( \begin{smallmatrix}
		-\frac{3452369}{4} & -\frac{25661603}{4}\\
		\frac{237917}{2} & \frac{1768447}{2}
		\end{smallmatrix}\right)$,\\
		
		 & $x^{-1}_2x_3x^{-3}_2,x^3_2x^{-1}_3x^{-1}_2x_3x^{-3}_2\}$ & 
		$\left( \begin{smallmatrix}
		\frac{80371}{4} & \frac{597401}{4}\\
		\frac{5145}{2} & \frac{38243}{2}
		\end{smallmatrix}\right)$,
		
		$\left( \begin{smallmatrix}
		\frac{1132425929}{4} & \frac{8417369243}{4}\\
		-\frac{152350279}{4} & -\frac{1132425989}{4}
		\end{smallmatrix}\right) \Biggr\}$\\
	& & \\
	\hline
	\end{tabular}
	\end{center}
	\end{small}
	\end{table}	
	
The Dealer obtains the sets
$$
U = \{u_1,u_2,u_3\} := \{ x_2x^{-1}_1x^{-1}_2x_3x^{-3}_2, x^{-1}_2x_3x^{-3}_2,
      x^3_2x_3^{-1}x^{-1}_2x_3x^{-3}_2\}	
$$
and
			\begin{align*}
			N&= \{N_1, N_2, N_3\}\\
			&:=\Biggl\{\begin{pmatrix}
				-\frac{3452369}{4} & -\frac{25661603}{4}\\
				\frac{237917}{2} & \frac{1768447}{2}
				\end{pmatrix},
				\begin{pmatrix}
				\frac{80371}{4} & \frac{597401}{4}\\
				\frac{5145}{2} & \frac{38243}{2}
				\end{pmatrix},
				\begin{pmatrix}
				\frac{1132425929}{4} & \frac{8417369243}{4}\\
				-\frac{152350279}{4} & -\frac{1132425989}{4}
				\end{pmatrix} \Biggr\}.
			\end{align*}
	
		\item He gets the shares 
			$(R_i,S_j)$ for the participants with $R_i \subset U$ and $S_j \subset N$ as follows:
 			\begin{enumerate}
			\item It is $m=\binom{n}{t-1}=\binom{3}{1}= 3$.
			\item The dealer chooses the elements $\tilde{a}_1, \tilde{a}_2, \tilde{a}_3$
			and gets the three sets 
			\begin{alignat*}{4}
			A_1 &= \left\{1\right\}, \qquad A_2 &= \left\{2\right\},\qquad A_3 &= \left\{3\right\}. 
			\end{alignat*}
			With the help of the $A_i$ the dealer gets the sets $R'_1, R'_2,$ and $R'_3$  which contain elements from the set $\left\{\tilde{a}_1, \tilde{a}_2, \tilde{a}_3\right\}$. He puts the element $\tilde{a}_{j}$ by which $i$ is not
			contained in the set $A_{j}$ for  
			$i=1,2,3$  
			and $j=1,2,3$, into the set $R'_{i}$.
			\begin{align*}
			1 \not\in A_2, A_3  \Longrightarrow  R'_1=\left\{\tilde{a}_2, \tilde{a}_3\right\},\\
			2 \not\in A_1, A_3 \Longrightarrow  R'_2=\left\{\tilde{a}_1, \tilde{a}_3\right\},\\
			3 \not\in A_1, A_2  \Longrightarrow  R'_3=\left\{\tilde{a}_1, \tilde{a}_2\right\}.
			\end{align*}
			
			Now we apply this to $U$ and $N$ to create the share-sets for the participants, respectively:  
			\begin{alignat*}{4}
			R_1 &= \left\{u_2, u_3\right\}, \qquad S_1 &= \left\{N_2, N_3\right\},\\
			R_2 &= \left\{u_1, u_3\right\}, \qquad S_2 &= \left\{N_1, N_3\right\},\\
			R_3 &= \left\{u_1, u_2\right\}, \qquad S_3 &= \left\{N_1, N_2\right\}.
			\end{alignat*}
			\end{enumerate}
\end{enumerate}	

\item The Dealer distributes to each participant  a tuple $(R_i,S_j)$. Participant $P_1$ gets $(R_1,S_2)$, $P_2$ gets $(R_2,S_3)$ and $P_3$ gets $(R_3,S_1)$.	
\end{enumerate}

Assume the participants $P_1$ and $P_2$ come together to reconstruct the secret. 
They generate the sets $U=\{u_1,u_2,u_3\}$ and $N=\{N_1,N_2,N_3\}$. 
The secret can be recovered as follows.

The participants apply regular Nielsen transformations step by step simultaneously 
for both sets $U$ and $N$ to get $X^{\pm}$ and $M^{\pm}$.  
The steps are shown in the Tables \ref{zuruckgehen1} and \ref{zuruckgehen2}.

\begin{table}[H]
\begin{small}
	\begin{center}
	\setlength{\tabcolsep}{.3\tabcolsep}
	\caption{Nielsen transformations (NTs) from the participants I}\label{zuruckgehen1}

	\begin{tabular}{|l | l| l|}
	\hline
	& & \\
	NTs & theoretical set U  & \hspace{1cm} explicit set N \\
	& & \\

	\hline
	\hline	
	& & \\
	 & $\{x_2x^{-1}_1x^{-1}_2x_3x^{-3}_2,$ &
	$\Biggl\{ \left( \begin{smallmatrix}
		-\frac{3452369}{4} & -\frac{25661603}{4}\\
		\frac{237917}{2} & \frac{1768447}{2}
		\end{smallmatrix}\right)$,\\
		
		 & $x^{-1}_2x_3x^{-3}_2,x^3_2x^{-1}_3x^{-1}_2x_3x^{-3}_2\}$ & 
		$\left( \begin{smallmatrix}
		\frac{80371}{4} & \frac{597401}{4}\\
		\frac{5145}{2} & \frac{38243}{2}
		\end{smallmatrix}\right)$,
		
		$\left( \begin{smallmatrix}
		\frac{1132425929}{4} & \frac{8417369243}{4}\\
		-\frac{152350279}{4} & -\frac{1132425989}{4}
		\end{smallmatrix}\right) \Biggr\}$\\
	& & \\

\hline	
	
		& & \\
	$(T1)_2$ & $\{x_2x^{-1}_1x^{-1}_2x_3x^{-3}_2,$ &
	$\Biggl\{ \left( \begin{smallmatrix}
		-\frac{3452369}{4} & -\frac{25661603}{4}\\
		\frac{237917}{2} & \frac{1768447}{2}
		\end{smallmatrix}\right) $,\\
		
		 & $x^{3}_2x^{-1}_3x_2,x^3_2x^{-1}_3x^{-1}_2x_3x^{-3}_2\}$ & 
		$\left( \begin{smallmatrix}
		\frac{38243}{2} & -\frac{597401}{4}\\
		-\frac{5145}{2} & \frac{80371}{4}
		\end{smallmatrix}\right)$,
		
		$\left( \begin{smallmatrix}
		\frac{1132425929}{4} & \frac{8417369243}{4}\\
		-\frac{152350279}{4} & -\frac{1132425989}{4}
		\end{smallmatrix}\right) \Biggr\}$\\
	& & \\
	\hline

	& & \\
	$(T2)_{3.2}$ & $\{x_2x^{-1}_1x^{-1}_2x_3x^{-3}_2,$ &
	$\Biggl\{ \left( \begin{smallmatrix}
		-\frac{3452369}{4} & -\frac{25661603}{4}\\
		\frac{237917}{2} & \frac{1768447}{2}
		\end{smallmatrix}\right)$,\\
		
		 & $\{x^{3}_2x^{-1}_3x_2,\, x^3_2x^{-1}_3\}$ & 
		 $\left( \begin{smallmatrix}
		\frac{38243}{2} & -\frac{597401}{4}\\
		-\frac{5145}{2} & \frac{80371}{4}
		\end{smallmatrix}\right)$,
		
		$\left( \begin{smallmatrix}
		5939 & 63664\\
		-799 & -8565
		\end{smallmatrix} \right) \Biggr\}$\\
	& & \\

	\hline

	& & \\
	$(T1)_{2}$ & $\{x_2x^{-1}_1x^{-1}_2x_3x^{-3}_2,$ &
	$\Biggl\{ \left( \begin{smallmatrix}
		-\frac{3452369}{4} & -\frac{25661603}{4}\\
		\frac{237917}{2} & \frac{1768447}{2}
		\end{smallmatrix}\right) $,\\
		
		 & $x^{-1}_2x_3x^{-3}_2,x^3_2x^{-1}_3\}$ & 
		$\left( \begin{smallmatrix}
		\frac{80371}{4} & \frac{597401}{4}\\
		\frac{5145}{2} & \frac{38243}{2}
		\end{smallmatrix}\right) $,
		
		$\left( \begin{smallmatrix}
		5939 & 63664\\
		-799 & -8565
		\end{smallmatrix}\right) \Biggr\}$\\
	& & \\
	\hline
	
	& & \\
	$(T2)_{2.3}$ & $\{x_2x^{-1}_1x^{-1}_2x_3x^{-3}_2,$ &
	$\Biggl\{ \left( \begin{smallmatrix}
		-\frac{3452369}{4} & -\frac{25661603}{4}\\
		\frac{237917}{2} & \frac{1768447}{2}
		\end{smallmatrix}\right) $,\\
		
		 & $x^{-1}_2,x^3_2x^{-1}_3\}$ & 
		 $\left( \begin{smallmatrix}
		-\frac{15}{2} & -\frac{221}{4}\\
		-1 & -\frac{15}{2}
		\end{smallmatrix}\right)$,
		
		$\left( \begin{smallmatrix}
		5939 & 63664\\
		-799 & -8565
		\end{smallmatrix}\right) \Biggr\}$\\
	& & \\
	\hline

	& & \\
	$(T2)_{1.3}$ & $\{x_2x^{-1}_1x^{-1}_2,x^{-1}_2,x^3_2x^{-1}_3\}$ &
	$\Biggl\{ \left( \begin{smallmatrix}
		\frac{653}{2} & \frac{9679}{4}\\
		-45 & -\frac{667}{2}
		\end{smallmatrix}\right) $,		
		  $\left( \begin{smallmatrix}
		-\frac{15}{2} & -\frac{221}{4}\\
		-1 & -\frac{15}{2}
		\end{smallmatrix}\right)$,
		$\left( \begin{smallmatrix}
		5939 & 63664\\
		-799 & -8565
		\end{smallmatrix}\right) \Biggr\}$\\
	& & \\
	\hline
	
	& & \\
	$(T1)_{2}$ & $\{x_2x^{-1}_1x^{-1}_2,x_2,x^3_2x^{-1}_3\}$ &
	$\Biggl\{ \left( \begin{smallmatrix}
		\frac{653}{2} & \frac{9679}{4}\\
		-45 & -\frac{667}{2}
		\end{smallmatrix}\right) $,		
		  $\left( \begin{smallmatrix}
		-\frac{15}{2} & \frac{221}{4}\\
		1 & -\frac{15}{2}
		\end{smallmatrix}\right)$,
		
		$\left( \begin{smallmatrix}
		5939 & 63664\\
		-799 & -8565
		\end{smallmatrix}\right) \Biggr\}$\\
	& & \\

		\hline
	\end{tabular}
	\end{center}
	\end{small}
	\end{table}

\begin{table}[H]
\begin{small}
	\begin{center}
	\setlength{\tabcolsep}{.3\tabcolsep}
	\caption{Nielsen transformations (NTs) from the participants II}\label{zuruckgehen2}

	\begin{tabular}{|l | l| l|}

\hline
	& & \\
	$(T2)_{1.2}$ & $\{x_2x^{-1}_1,x_2,x^3_2x^{-1}_3\}$ &
	$\Biggl\{ \left( \begin{smallmatrix}
		-29 & -109\\
		4 & 15
		\end{smallmatrix}\right)$,		
		  $\left( \begin{smallmatrix}
		-\frac{15}{2} & \frac{221}{4}\\
		1 & -\frac{15}{2}
		\end{smallmatrix}\right)$,
		
		$\left( \begin{smallmatrix}
		5939 & 63664\\
		-799 & -8565
		\end{smallmatrix}\right) \Biggr\}$\\
	& & \\
	\hline
	
	& & \\
	$(T1)_{1}$ & $\{x_1x^{-1}_2,x_2,x^3_2x^{-1}_3\}$ &
	$\Biggl\{ \left( \begin{smallmatrix}
		15 & 109\\
		-4 & -29
		\end{smallmatrix}\right) $,		
		  $\left( \begin{smallmatrix}
		-\frac{15}{2} & \frac{221}{4}\\
		1 & -\frac{15}{2}
		\end{smallmatrix}\right)$,
		
		$\left( \begin{smallmatrix}
		5939 & 63664\\
		-799 & -8565
		\end{smallmatrix}\right) \Biggr\}$\\
	& & \\
	\hline
	
	& & \\
	$(T2)_{1.2}$ & $\{x_1,x_2,x^3_2x^{-1}_3\}$ &
	$\Biggl\{ \left( \begin{smallmatrix}
		-\frac{7}{2} & \frac{45}{4}\\
		1 & -\frac{7}{2}
		\end{smallmatrix}\right)$,		
		  $\left( \begin{smallmatrix}
		-\frac{15}{2} & \frac{221}{4}\\
		1 & -\frac{15}{2}
		\end{smallmatrix}\right)$,
		
		$\left( \begin{smallmatrix}
		5939 & 63664\\
		-799 & -8565
		\end{smallmatrix}\right) \Biggr\}$\\
	& & \\
	\hline

	& & \\
	$(T1)_{3}$ & $\{x_1,x_2,x_3x^{-3}_2\}$ &
	$\Biggl\{ \left( \begin{smallmatrix}
		-\frac{7}{2} & \frac{45}{4}\\
		1 & -\frac{7}{2}
		\end{smallmatrix} \right)$,		
		  $\left( \begin{smallmatrix}
		-\frac{15}{2} & \frac{221}{4}\\
		1 & -\frac{15}{2}
		\end{smallmatrix}\right)$,
		
		$\left(\begin{smallmatrix}
		-8565 & -63664\\
		799 & 5939
		\end{smallmatrix}\right) \Biggr\}$\\
	& & \\
	\hline
	
	& & \\
	$\lbrack (T2)_{3.2} \rbrack^3$ & $\{x_1,x_2,x_3\}$ &
	$\Biggl\{ \left( \begin{smallmatrix}
		-\frac{7}{2} & \frac{45}{4}\\
		1 & -\frac{7}{2}
		\end{smallmatrix}\right)$,		
		  $\left( \begin{smallmatrix}
		-\frac{15}{2} & \frac{221}{4}\\
		1 & -\frac{15}{2}
		\end{smallmatrix}\right) $,
		
		$\left( \begin{smallmatrix}
		-11 & 120\\
		1 & -11
		\end{smallmatrix}\right) \Biggr\}$\\
	& & \\
	\hline
	\end{tabular}
	\end{center}
	\end{small}
	\end{table}

 With the knowledge of the set $M^{\pm}= \Biggl\{ \left( \begin{smallmatrix}
		-\frac{7}{2} & \frac{45}{4}\\
		1 & -\frac{7}{2}
		\end{smallmatrix}\right)$,		
		  $\left( \begin{smallmatrix}
		-\frac{15}{2} & \frac{221}{4}\\
		1 & -\frac{15}{2}
		\end{smallmatrix}\right) $,
		$\left( \begin{smallmatrix}
		-11 & 120\\
		1 & -11
		\end{smallmatrix}\right) \Biggr\}$ the participants can reconstruct the secret easily.
It is  
	\begin{align*}
	a_1:=tr(M_1)=-7, \quad a_2:=tr(M_2)=-15, \quad a_3:=tr(M_3)=-22
	\end{align*}
	and hence it is
	\begin{align*}
	S:= \sum_{j=1}^{3} \frac{1}{|a_j|} 
    =  \frac{1}{7} + \frac{1}{15} + \frac{1}{22} 
	= \frac{589}{2310}.
	\end{align*}

\end{example}

In general we can use any free matrix group~$F$ of rank $m:=\binom{n}{t-1}$ for a $(n,t)$-secret sharing scheme as it is described in this section.
The shares can be generated by the above method and are tuples $(R_i,S_j)$ with 
$R_i \subset U$ and $S_j \subset N$. Some other ideas for the secret~$S$ are

\begin{align*}
S&:=\prod_{i=1}^{m}|tr(M_i)| \text{ or } S:=  \sum_{i=1}^{m}|tr(M_i)| \text{ or }\\
S&:= \prod_{i=1}^{m}(tr(M_i))^2 \text{ or } S:=  \sum_{i=1}^{m} (tr( M_i))^2 \text{ or } \\
S&:=  \prod_{i=1}^{\frac{m}{2}}tr(\lbrack M_{2i-1},M_{2i}\rbrack) \text{ if $m$ is even} \text{ or }
S:= \sum_{i=1}^{m} tr(M_i^2).
\end{align*}


\section{A variation of the secret sharing scheme based on Nielsen transformations}%
\label{varsecretsharingNT}

We explain a variation of the secret sharing scheme described in Section \ref{NSec}.
As in the previous sections, let $F$ be a free group with the abstract free generating set $X:=\{x_1,x_2,\ldots,x_q\}$, $q \in \mathbb N$, that is,
\begin{align*}
F= \langle X ; \phantom {R} \rangle.
\end{align*}

For a $(n,t)$-secret sharing scheme the dealer chooses a Nielsen reduced set $U:=\{u_1,u_2,\ldots,u_m\}  \subset F$, with $m=\binom{n}{t-1}$. The $u_i$ are given as words in $X$.
The secret is the sum 
\begin{align*}
S:=\sum_{i=1}^{m} \frac{1}{|u_i|},
\end{align*}
with $| u_i |$ the length of the word $u_i$.

The dealer does a regular Nielsen transformation on the set $U$ to get the Nielsen-equivalent set $V$ (Figure \ref{regNielsentanfs1}).  

\begin{figure}[H]
			\begin{center}			
			\begin{tikzpicture}
			\node (a) at (0.5,1.5) {$U:=\{ u_1, u_2,\ldots, u_m\}$}; 
			\draw[->] (0,1)--(0,-1);
			\node (s) at (1.5,0) {regular Nielsen};
			\node (t) at (1.5,-0.5) {transformation};
			\node (b) at (0.5,-1.5) {$V:=\{ v_1, v_2,\ldots, v_m\}$};
			\end{tikzpicture}
			\caption{Regular Nielsen transformation}\label{regNielsentanfs1}
			\end{center}
\end{figure}

Each participant $P_i$ , $1 \leq i \leq n$, gets one set $R_i \subset V$, as above.

If $t$ of the $n$ participants come together to reconstruct the secret, they combine their shares and get the set $V=\{v_1,v_2,\ldots,v_m\}$. They have to find a Nielsen-reduced set $U':=\{u'_1,u'_2,\ldots,u'_m\}$ to $V$. They apply Nielsen transformations in a Nielsen reducing manner as described in \cite{CRR} and \cite{LS} and get from $V$ a Nielsen-reduced set $U'$. The secret is the sum
\begin{align*}
S=\sum_{i=1}^{m} \frac{1}{|u'_i|},
\end{align*}
because $\sum_{i=1}^{m} |u'_i|=\sum_{i=1}^{m} |u_i|$ for the Nielsen reduced sets $U'$ and $U$ (see \cite[Corollary 2.9]{CRR}).


\section{A symmetric key cryptosystem using Nielsen transformations}\label{cryptoNielsen}

Before Alice and Bob are able to communicate with each other they have to make some arrangements.
Let $F$ be an abstract free group with the free generating set \newline $X=\{x_1,x_2,\ldots,x_q\}$,  $q \in \mathbb N  \backslash \{1\} $. Let
\begin{align*}
 \varphi: F &\rightarrow SL(2, \mathbb Q)\\
 x_i &\mapsto M_i,
\end{align*}
be a faithful representation of $F$ into $SL(2,\mathbb Q)$ as in Section \ref{NSec}. The group $G=\varphi(F)$ is isomorphic to $F$ under the map  $x_i \mapsto M_i$, for $i=1,\ldots,q$.\\
Let $N$ be the number of letters from the alphabet $A=\{a_1,\ldots,a_N\}$, for instance $N=26$. We assume that $N\geq 5$.\\
Let $U\subset F$, $U=\{u_1,\ldots,u_N\}$ be a basis of a free subgroup of $F$ of rank $N$. Such systems $U$ are easily to construct (see \cite{CRR} or \cite{LS}).\\
There is the one to one assignment 
\begin{align*}
 A &\rightarrow U\\
 a_{j} &\mapsto u_{j}, \quad \text{ for } j=1,\ldots,N.
\end{align*}
Let $U'=\varphi(U)=\{U'_1,\ldots,U'_N\} \subset SL(2,\mathbb Q)$, $u_{j} \mapsto U_{j}'$ for $j=1,\ldots,N$.
The set $U'$ is a basis for a free subgroup of $G$. Now, Alice and Bob agree on a block sequence $P:=p_1p_2\cdots p_k$ with, say, $1 \leq p_i \leq 4$ and $k \geq 2$, and for each $p_i$ they construct a regular Nielsen transformation $f_i$ from $U'$ to a Nielsen-equivalent system $f_i(U')=\{V'_{i_1}, \ldots, V'_{i_N}\}, f_i(U'_{j})=V_{i_j}'$, $j=1,\ldots,N$.
 The Nielsen transformations $f_i$, $1 \leq i \leq k$, should be pairwise different and given as sequences of elementary Nielsen transformations from $U'$ to $f_i(U')$.
 
\vspace{0.5cm}

As soon as Alice an Bob agree on a Nielsen transformation $f_i$ they compute $f_i(U')$, $i=1,2,\ldots,k$, independent from each other even if they do not know the message. 
Hence they get a one to one assignment between the letters in their alphabet $A$ and the matrices for the ciphertext depending from the part of the sequence $P$. This is shown in Table~\ref{ZuordnungCST}.

\begin{table}[htb]
\caption{Assignment: Letters in $A=\{a_1,a_2,\ldots,a_N\}$ to matrices for the ciphertext depending from the part of the Sequence $P$}\label{ZuordnungCST}
\[
\begin{array}{c c| c c c c c }
& &  \multicolumn{5}{c}{\text{Elements from the alphabet $A$}}\\
& & &  & & &   \\
& & & a_1 & a_2 & \ldots & a_{N}\\
& & &  & & &  \\
\hline
& & &  & & &  \\
& \text{Sequence P} & & & & &  \\
& & &  & & &  \\
&  p_1 & & V'_{1_1} & V'_{1_2} & \ldots & V'_{1_N}  \\
& & &  & & &  \\
&  p_2 &  & V'_{2_1} & V'_{2_2} & \ldots & V'_{2_N} \\
& & & & & &  \\
& \vdots &  & \vdots & \vdots & \ldots & \vdots \\
& & & & & & \\
&  p_k &  & V'_{k_1} & V'_{k_2} & \ldots & V'_{k_N} \\
& & & & & & \\

\end{array}
\]
\end{table}

\vspace{0.5cm}

 Now, Alice wants to send a message $S$ with $z$, $z>0$, letters from $A$. 
 To describe the procedure let first 
 \begin{align*}
 z=|S|=\sum_{i=1}^{k}p_i.
 \end{align*}
 Alice cuts the message $S$ into pieces corresponding to the sequence $P$, that is,
 \begin{align*}
 S= a_{1_{1}}\cdots a_{{1}_{p_1}} a_{2_{1}}\cdots a_{{2}_{p_2}}  \cdots  a_{k_{1}}\cdots a_{{k}_{p_k}},
 \end{align*}
 with $a_{i_j} \in A$ for $1 \leq i \leq k$ and $1 \leq j \leq p_i$.
 Then she uses the Table \ref{ZuordnungCST} to get the matrices for the ciphertext depending on the sequence $p_i$, $1 \leq i \leq k$:
\begin{align*}
\text{If } a_{i_j} = a_{t} \qquad \text{then } a_{i_j} \mapsto V'_{i_t},  \quad 1\leq t \leq N, \quad 1 \leq j \leq p_i.
\end{align*}
The ciphertext $C$ is the following sequence of matrices:
\begin{align*}
 C=V_{1_1}\cdots V_{1_{p_1}}  V_{2_1}\cdots V_{2_{p_2}} \cdots  V_{k_1}\cdots V_{k_{p_k}},
 \end{align*}
 with $V_{g_q} \in \{V'_{g_1},V'_{g_2}, \ldots, V'_{g_N}\}$ for $1 \leq g \leq k$ and $1 \leq q \leq p_g$. 
 Alice sends the ciphertext $C$ just as a sequence of matrices to Bob. To reconstruct the message $S$ Bob does the following steps:
\begin{enumerate}
\item He cuts back the ciphertext into pieces as above corresponding to the known sequence $P$.
\item Because he gets the same table as Alice (Table \ref{ZuordnungCST}) he can match $V'_{i_j}$ to $a_t$  for the piece corresponding to  $p_i$, with $1\leq i \leq k$, $1 \leq j\leq p_i$ and $1 \leq t \leq N$, and hence reads the message $S$ in the alphabet easily.
\end{enumerate}

Now, assume that $z=|S|\ne \sum_{i=1}^{k} p_i$.

If $z < \sum_{i=1}^{k} p_i$, there is no problem, Alice just ends with the last letter. 

If $S=\tilde{S}_1\tilde{S}_2$ with $|\tilde{S}_1|=\sum_{i=1}^{k}p_i$ then she first applies the above procedure to $\tilde{S}_1$ and continues then with $\tilde{S}_2$ in the same manner.

Indeed, if $z=|S| > \sum_{i=1}^{k}p_i$ and $S=\tilde{S}_1\tilde{S}_2 \cdots \tilde{S}_m$ with $|\tilde{S}_1|= \cdots =|\tilde{S}_m|$ then we may improve the cryptosystem. Alice and Bob agree in addition on a permutation $\sigma \in S_m$, $S_m$ the symmetric group on $m$ symbols, and work with $\tilde{S}_{\sigma(1)}\tilde{S}_{\sigma(2)} \cdots \tilde{S}_{\sigma(m)}$ instead of $S=\tilde{S_1}\tilde{S_2} \cdots \tilde{S}_{m}$. Bob starts then the decryption procedure with applying first $\sigma^{-1}$.

\begin{remark}
\begin{enumerate}
\item If Alice wants to send several messages to Bob or vice versa, then to improve the system, they may replace during each message transmission the Nielsen transformations $f_i$ by different Nielsen transformations, for instance by $\tilde{f}_i=f_if^k$ for a fixed Nielsen transformation $f$ on $U$ and $k=0,1,2,\ldots$. They also may replace the system $U \subset F$ by different systems, such that $U$ is not used too often.
\item The cryptosystem is a polyalphabetic system. A matrix $u \in U$, and hence a letter $a \in A$, can be encrypted differently at different times.
\item A possible attacker Eve cannot read the message. She can only see a sequence of matrices in $SL(2,\mathbb Q)$. From this she neither knows the system $U'=\varphi(U)$ nor the Nielsen transformations $f_i$.
Since the block lengths $p_i$ are very small (it is $ 1 \leq p_i \leq 4$ for all $i$) and since rather frequently the Nielsen transformations $f_i$ and the system $U$ will be changed, a statistical frequency attack is almost impossible.\\
Moreover, if we just have two $N$-tuples $\{A_1,\ldots,A_N\}$ and $\{B_1,\ldots,B_N\}$ in $SL(2,\mathbb Q)$, it is hard to decide if they generate the same group, in fact, if $A \in SL(2, \mathbb Q)$ and $H$ a subgroup of $SL(2, \mathbb Q)$, it is hard to decide if $A \in H$.
\end{enumerate}
\end{remark}

More details and generalizations of the cryptosystem  as well as 
more cryptographical analysis can be found in \cite{M} by A. Moldenhauer.

\section{Cryptosystem with Nielsen transformation inspired by the ElGamal cryptosystem}\label{ElGamal}
Now we describe a public key cryptosystem for Alice and Bob which is inspired by the ElGamal cryptosystem (see \cite{E} or \cite[Section 1.3]{MSU}), based on discrete logarithms, that is:
\begin{enumerate}
\item Alice and Bob agree on a finite cyclic group $G$ and a generating element $g \in G$.
\item Alice picks a random natural number $a$ and publishes the element $c:=g^a$.
\item Bob, who wants to send a message $m \in G$ to Alice, picks a random natural number $b$ and sends the two elements $m\cdot c^b$ and $g^b$, to Alice. Note that $c^b=g^{ab}$.
\item Alice recovers $m=\left( m \cdot c^b\right) \cdot \left( \left( g^b \right)^a\right)^{-1}$. 
\end{enumerate}

In the following cryptosystem let $N \in \mathbb N$ be the number of  letters from the used alphabet, $X=\{x_1,x_2\ldots,x_N\}$, $S=X \cup X^{-1}$ and $F$ the free group on the free generating set $X$, that is,  $F=\langle X ; \phantom{R} \rangle$.  The message is an element $m \in S^*$, the set of all words on $S$.
Public are the free group $F$, its generating set $X$, an element $a \in S^*$
and an automorphism $f: F \rightarrow F$, of infinite order, given as a Nielsen transformation or a Whitehead automorphism (see \cite{DKR}).
Each automorphism of $F$ is a product of elementary Nielsen transformations between two bases of $F$ (see \cite[Korollar 2.10]{CRR}).\\

The cryptosystem is now as follows:\\

\underline{Public parameters:} The group $F= \langle X ; \phantom{R} \rangle$, an element $a \in F$ and an automorphism $f: F \to F$ of infinite order.

\begin{enumerate}
\item Alice chooses a $n \in \mathbb N$ and publishes the element $c:=f^n(a) \in S^*$.

\item Bob picks a random $t \in \mathbb N$ and his message $m\in S^*$.
He calculates $c_1:=m \cdot f^t(c) \in S^*$ and $c_2:=f^t(a) \in S^*$. He sends  the ciphertext $(c_1,c_2) \in S^* \times S^*$ to Alice. 

\item Alice calculates
\begin{align*}
c_1 \cdot f^n(c_2)^{-1}&= m\cdot f^t(c) \cdot f^n(c_2)^{-1}\\
                          &= m\cdot f^t(f^n(a)) \cdot (f^n(f^t(a))^{-1}\\
                          &= m\cdot f^{t+n}(a) \cdot (f^{n+t}(a))^{-1}\\
                          &= m,
\end{align*}
and gets the message $m$.
\end{enumerate}

\begin{remark}
A possible attacker, Eve,  can see the elements  $c, c_2, c_1 \in S^*$. She does not know the free length of $m$ and the cancellations between $m$ and $f^t(c)$ in  $c_1$. Hence she gets no information about $m$ from the element $c_1$.   Eve just sees words in the free generating system from which it is almost impossible to realize the exponents $n$ and $t$, that is, the private keys from Alice and Bob, respectively.
\end{remark}

\begin{remark} We give some ideas to  enhance the security, they can also be combined:
\begin{enumerate}
\item The element $a \in F$ is a common secret between Alice and Bob. They could use for example the Anshel-Anshel-Goldfeld key exchange protocol (see \cite{MSU}) to agree on the element $a$.
\item Alice and Bob agree on a faithful representation from $F$ into the special linear group of all $2 \times 2$ matrices with entries in $\mathbb Q$, that is, $g: F \to SL(2, \mathbb Q)$. Now $m \in S$ and Bob sends $c_1:=g(m)  \cdot g(f^t(c)) \in SL(2, \mathbb Q)$ instead of $c_1:=m \cdot f^t(c) \in S^*$. Therefore Alice calculates $c_1 \cdot (g(f^n(c_2)))^{-1}=g(m)$ and hence the message $m=g^{-1}(g(m))\in S$.
\end{enumerate}    
\end{remark}



\bibliographystyle{abbrv} 

\section*{Author information}
Anja I. S. Moldenhauer, Fachbereich Mathematik, Universität Hamburg, Bundesstrasse 55,\\
20146 Hamburg, Germany.\\
E-mail: \texttt{anja.moldenhauer@uni-hamburg.de}\\

Gerhard Rosenberger, Fachbereich Mathematik, Universität Hamburg, Bundesstrasse 55,\\
20146 Hamburg, Germany.\\
E-mail: \texttt{gerhard.rosenberger@math.uni-hamburg.de}\\

\end{document}